\documentclass{amsart}

\usepackage{stmaryrd}

\usepackage{amsmath, amsthm, amssymb}

\theoremstyle{definition}

\usepackage{scalerel}
\usepackage{mathrsfs}

\usepackage{tikz}
\usepackage{hyperref}
\usepackage{placeins}

\usepackage{proof}

\usepackage{verbatim}
\usepackage{mathtools}
\usepackage{cleveref}

\usepackage{bbm}

\usepackage[bottom]{footmisc}

\usepackage{enumitem}

\usepackage{biblatex} 
\newtheorem{theorem}{Theorem}
\newtheorem*{theorem*}{Theorem}
\newtheorem*{maintheorem*}{Main Theorem}

\newtheorem{definition}[theorem]{Definition}
\newtheorem{corollary}{Corollary}[theorem]
\newtheorem{prop}[theorem]{Proposition}
\newtheorem{lemma}[theorem]{Lemma}

\newtheorem{claim}{Claim}[theorem]

\newenvironment{claimproof}[1]{\par\noindent\textit{Proof of the Claim.}\space#1}{\hfill $\blacksquare$}

\newtheorem{manualtheoreminner}{Theorem}
\newenvironment{manualtheorem}[1]{%
  \IfBlankTF{#1}
    {\renewcommand{\themanualtheoreminner}{\unskip}}
    {\renewcommand\themanualtheoreminner{#1}}%
  \manualtheoreminner
}{\endmanualtheoreminner}

\newcommand\A{\mathscr{A}}
\newcommand\I{\mathscr{I}}
\newcommand\G{\mathcal{G}}
\newcommand\F{\mathcal{F}}

\newcommand\carda{\mathfrak{a}}
\newcommand\cardb{\mathfrak{b}}
\newcommand\cardc{\mathfrak{c}}
\newcommand\cardd{\mathfrak{d}}

\newcommand\Q{\mathbb{Q}}

\newcommand\ZFU{\textup{ZFU}}
\newcommand\ZFCU{\textup{ZFCU}}
\newcommand\ZFUR{\textup{ZFU}_\textup{R}}

\def\<#1>{\left\langle#1\right\rangle}
\def\[#1]{\left\llbracket#1\right\rrbracket}

\title{\textbf{Plenitudinous Urelements and the Definability of Cardinality}}
\author{Bokai Yao}\address[Bokai Yao]{Peking University}
 \email{bkyao@pku.edu.cn}
\urladdr{https://bokaiyao.com}

\thanks{Part of this paper appeared in the author’s doctoral dissertation at the University of Notre Dame. The author thanks Joel David Hamkins for helpful discussions and Elliot Glazer for identifying and fixing an error in an earlier draft.}

\usepackage{graphicx}

\begin{document}
\maketitle

\begin{abstract}
The Axiom of Plenitude asserts that every ordinal is equinumerous to a set of urelements, while its stronger form, Plenitude$^+$, extends it to all sets. We investigate these two axioms within ZF set theory with urelements. Assuming that cardinality is definable, Plenitude$^+$ together with the Collection Principle implies the Reflection Principle. If either cardinality is representable or Small Violations of Choice (SVC) holds, Plenitude$^+$ implies the Reflection Principle. In contrast, Plenitude is considerably weaker: SVC + Plenitude does not prove the Collection Principle, and SVC + Plenitude + Reflection Principle does not prove Plenitude$^+$.
\end{abstract}

\section{Introduction}

Maximality principles in set theory often assert, in one way or another, that there are as many sets as possible. If the set-theoretic universe also includes \textit{urelements}---elements of sets that are not sets---then it is natural to consider maximality principles concerning the number of urelements. One natural principle of this kind asserts that there are as many urelements as ordinals. That is,
\begin{itemize}
\item [] (Plenitude) For every ordinal $\alpha$, there is a set of urelements equinumerous to $\alpha$.
\end{itemize}
In the absence of the Axiom of Choice, a natural strengthening of Plenitude is to assert that there are as many urelements as sets.
\begin{itemize}
\item [] (Plenitude$^+$) For every set $x$, there is a set of urelements equinumerous to $x$.
\end{itemize}
In this paper, we study these two principles in ZFU, i.e., ZF set theory with urelements. In particular, we investigate their relationship with standard set-theoretic axioms such as the  Collection Principle and the Reflection Principle. We obtain positive results using a weak but natural choice principle in urelement set theory: the definability of cardinality. The main positive results, proved in Section \ref{section:ZFUResults}, are the following.
\begin{theorem}
Assume that cardinality is parametrically definable and Plenitude$^+$ holds. The Collection Principle is equivalent to the Reflection Principle.
\end{theorem}

\begin{theorem}
Assume that cardinality is parametrically representable. Plenitude$^+$ implies the Reflection Principle.
\end{theorem}
\noindent Furthermore, we consider the principle of \textit{Small Violations of Choice} (SVC), a central principle introduced by Blass \cite{blassInjec1979}, which asserts that the Axiom of Choice holds in some forcing extension.
\begin{theorem}
Assume SVC. Plenitude$^+$ implies the Reflection Principle.
\end{theorem}
\noindent We prove the following independent results in Section \ref{section:IndependenceResults}, showing that Plenitude is considerably weaker in the choiceless context.
\begin{theorem}\label{thm:Independence1}
$\ZFU$ + Plenitude + SVC does not prove the Collection Principle.
\end{theorem}

\begin{theorem}\label{thm:Independence2}
$\ZFU$ + Plenitude + SVC + Reflection Principle does not prove Plenitude$^+$.
\end{theorem}
\noindent Using the model constructed in the last theorem, we show that SVC, although it implies that cardinality is parametrically definable, does not imply the parameter-free definability of cardinality.

\section{Preliminaries}
In this section, we introduce our notation, review some key facts and known results.

The language of urelement set theory, in addition to $\in$, contains a unary predicate $\A$ for urelements. The axioms of $\ZFU$\footnote{In \cite{yao2023set}, \cite{YaoAxiomandForcing} and \cite{YaoAbstraction}, this theory is called  $\ZFUR$, where R indicates that the theory only includes Replacement and not Collection. However, most authors assume by default that the standard axioms of ZF with urelements include Replacement but not Collection, which makes the subscript R unnecessary and somewhat confusing. Therefore, we follow the standard notation and no longer use $\ZFUR$.} include Foundation, Pairing, Union, Powerset, Infinity, Separation, Replacement, Extensionality \textit{for sets}, and the axiom that no urelements have members (see \cite[Section 1]{YaoAxiomandForcing} for the precise formulation of these axioms). Importantly, $\ZFU$ allows a proper class of urelements. $\ZFCU$ is $\ZFU$ + the Axiom of Choice, while ZF(C) is ZF(C)U + $\forall x \neg \A(x)$.

In $\ZFU$, every object $x$ has a \textit{kernel}, denoted by $ker(x)$, which is the set of the urelements in $TC(\{x\})$ (the transitive closure of $\{x\}$). A set is pure if its kernel is empty. $V$ denotes the class of all pure sets. $Ord$ is the class of all ordinals, which are transitive \textit{pure} sets well-ordered by the membership relation. For any sets $x$ and $y$, $x \sim y$ ($x$ is equinumerous to $y$) if there is a bijection between $x$ and $y$; $x \preceq y$ if there is an injection from $x$ into $y$; and $x \preceq^* y$ if either there is a surjection from $y$ onto $x$ or $x = \emptyset$. An ordinal is \textit{initial} if it is not equinumerous to any smaller ordinals. For every set $x$, there is the least initial ordinal that does not inject into $x$, which is known as \textit{the Hartogs number} and denoted by $\aleph(x)$. $\A$ will also stand for the class of all urelements. $A \subseteq \A$ thus means ``$A$ is a set of urelements''. For any $A \subseteq \A$, the $V_{\alpha}(A)$- hierarchy is defined in the standard way, i.e., 
\vspace{4pt}
\begin{itemize}
    \item [] $V_0(A) = A$;
    \item [] $V_{\alpha+1}(A) = P(V_{\alpha}(A)) \cup A$;
    \item [] $V_{\gamma}(A) = \bigcup_{\alpha < \gamma} V_\alpha(A)$, where $\gamma$ is a limit;
    \item [] $V(A) = \bigcup_{\alpha \in Ord} V_\alpha(A)$.
\end{itemize}
\vspace{4pt}
For every object $x$ and set $A \subseteq \A$, $x \in V(A)$ if and only if $ker(x) \subseteq A$. $U$ denotes the class of all objects, i.e., $U = \bigcup_{A\subseteq\A} V(A)$.

An important feature of $\ZFU$ is that $U$ has many definable nontrivial automorphisms. In particular, every permutation $\pi$ of a set of urelements can be extended to a definable permutation of $\A$ by letting $\pi$ be identity elsewhere, and then $\pi$ can be further extended to a permutation of $U$ by letting $\pi x$ be $\{\pi y : y \in x \}$ for every set $x$. $\pi$ preserves $\in$ so it is an automorphism of $U$. For every $x$ and automorphism $\pi$ of $U$, $\pi$ pointwise fixes $x$ whenever $\pi$ pointwise fixes $ker(x)$.

Although $\ZFCU$ looks very much like \textit{the} urelement analog of ZFC, there turns out to be a hierarchy of axioms that are independent of $\ZFCU$ (\cite[Theorem 2.5]{YaoAxiomandForcing}), some of which are ZF-theorems. In particular, the following two are not provable from $\ZFCU$.\\

\begin{itemize}
\item[] (Collection) $\forall w, u (\forall x \in w \ \exists y \varphi(x, y, u)   \rightarrow \exists v \forall x \in w\  \exists y \in v\  \varphi(x, y, u))$.\\
\item [] (RP) For every set $x$ there is a transitive set $t$ with $x \subseteq t$ such that for every $x_1, ..., x_n \in t$, $\varphi(x_1, ..., x_n) \leftrightarrow \varphi^t(x_1, ..., x_n)$.\\
\end{itemize}

\noindent The following observation will be useful later on.
\begin{prop}\label{prop:CollectionEquivs}
The following are equivalent over $\ZFU$
\begin{enumerate}
\item Collection
\item $\forall w, u (\forall x \in w \ \exists y \varphi(x, y, u)   \rightarrow \exists A \subseteq \A \forall x \in w\  \exists y \in V(A) \  \varphi(x, y, u))$
\end{enumerate}
\end{prop}
\begin{proof}
It is clear that (1) $\to$ (2) since every set is in some $V(A)$. Assume (2) and suppose that $\forall x \in w \ \exists y \varphi(x, y, u)$. Then there is some set of urelements $A$ such that $\forall x \in w\  \exists y \in V(A) \  \varphi(x, y, u)$. For every $x \in w$, let $\alpha_x$ be the least ordinal such that $\exists y \in V_\alpha(A) \varphi(x, y, u)$ and let $\alpha = \sup_{x \in w} \alpha_x$. Then $V_\alpha (A)$ is the desired collection set.
\end{proof}
\noindent Thus, according to Collection, the universe always contains a set of urelements large enough to generate the required collection set. Whereas Replacement is often viewed as a height-maximality principle about the universe, Collection can be understood as a maximality principle asserting that the universe is both tall and wide.

The idea of reflection, which asserts that the universe is indescribable, underlies many maximality principles in set theory, including certain large cardinal axioms. The reflection principle RP formulated above is a natural modification of the L\'evy-Montague Principle, which is a fundamental theorem of ZF. RP implies Collection over $\ZFU$ by a standard argument (\cite[Proposition 11]{yao2023set}). Moreover, it is shown in \cite[Theorem 2.5]{YaoAxiomandForcing} that Collection is equivalent to RP over $\ZFCU$, and it is proved in \cite{GlazerYaoRPinZFU} that Collection does not imply RP over ZFU.

On the other hand, over $\ZFCU$ enough urelements yield RP.
\begin{definition}
A set $x$ is \textit{realized by} $A$ if  $x \sim A$; a set is \textit{realized} if it is realized by some set of urelements.  \textit{Plenitude} is the principle that every ordinal is realized. \textit{Plenitude}$^+$ is the principle that every set is realized.
\end{definition}
\noindent We note that $\ZFCU$ + Plenitude is equiconsistent with ZF (\cite[Theorem 10]{yao2023set}) and that over $\ZFCU$, ``the urelements form a proper class'' does not imply Plenitude. In fact, there are models of $\ZFCU$ + RP in which the urelements form a proper class but every set of urelements is countable (see \cite[Theorem 2.17 (1)]{YaoAxiomandForcing}). It is proved in \cite[Theorem 2.5]{YaoAxiomandForcing} that Plenitude implies RP (and hence Collection) over $\ZFCU$. In other words, with AC the maximality principles about urelements imply certain maximality principles about $U$. The use of AC is necessary: it is proved in \cite{GlazerYaoRPinZFU} that Plenitude$^+$ does not imply RP over ZFU.

Given these known results, one may ask if there are weaker but natural choice principles that can rebuild the connections among these maximality principles. We now turn to one such principle---the definability of cardinality.

\begin{definition}\
\begin{enumerate}
\item \textit{Cardinality is definable} if some formula $\varphi(x, |x|)$ defines a \textit{cardinality assignment (class) function} $x \mapsto |x|$ such that for every sets $x$ and $y$, $x \sim y$ if and only if $|x| = |y|$. 

\item Cardinality is \textit{parametrically} definable if some formula $\varphi(x, |x|, p)$ together with some object $p$ defines a cardinality assignment function.
\end{enumerate}
\end{definition}
\noindent Given a cardinality assignment function $x \mapsto |x|$, every object in the range of this function is said to be a \textit{cardinal}, denoted by Frankfurt font letters $\carda, \cardb, \cardc, \cardd ...$ Cardinal addition $\carda + \cardb$ and $\carda \leq \cardb$ are defined as usual.

Cardinality is always definable if either the urelements form a set or AC holds. It is consistent with $\ZFU$ that cardinality is not parametrically definable, as proved independently by Gauntt \cite{gauntt1967undefinability} and L\'evy \cite{levy1969definability}. L\'evy \cite{levy1969definability} also shows that it is consistent with $\ZFU$ that cardinality is parametrically definable but not definable. The definability of cardinality thus becomes a choice principle in urelement set theory. We conclude this section with a useful lemma.

\begin{lemma}\label{lemma:HP->AlmostUniversalSet}\
\begin{enumerate}
\item Suppose that cardinality is parametrically definable. For every cardinal $\carda$, $\{ \cardb \mid \cardb \leq \carda\}$ is a set.

\item If cardinality is parametrically definable, then there is a set $A$ of urelements such that every cardinal is in $V(A)$;

\item  If cardinality is definable, then either $\A$ is a set or every cardinal is in $V$.
\end{enumerate}
\end{lemma}
\begin{proof}
(1) Let $x$ be such that $|x| = \carda$. If $\cardb \leq \carda$, then $|z| = \cardb$ for some $z \preceq x$ and so $|z| = |y|$ for some $y \subseteq x$. It follows that $\{ \cardb \mid \cardb \leq \carda\} = \{ |y| \mid y \subseteq x \}$, and the latter is a set since it is the range of the cardinality assignment function on $P(x)$.

(2) We may assume that $\A$ is a proper class. Suppose that some formula $\varphi(x, |x|, p)$, together with an object $p$, defines a cardinality assignment function. Let $A = ker(u)$. Suppose that $|x| \notin V(A)$ for some set $x$. Fix some $a \in ker(|x|) - A$ and some $b \notin A \cup ker(|x|)$. Let $\pi$ be the automorphism that only swaps $a$ and $b$. So $\pi |x| \neq |x|$. We also have $\varphi(\pi x, \pi |x|, u)$ so $\pi |x| = | \pi x|$. But $\pi x \sim x$, so $ |\pi x| = |x|$---contradiction.

(3) Suppose that cardinality is definable and $\A$ is a proper class. Then it follows that every cardinal is in $V$ by letting $A$ be $\emptyset$ as in (2). \end{proof}

\section{$\ZFU$ Results}\label{section:ZFUResults}

\subsection{Duplication}
We now work in $\ZFU$. Our key lemma will use the assumption that we may \textit{duplicate} urelements whenever needed, which motivates the following definition.
\begin{definition}
Let $A$ be a set of urelements.
\begin{enumerate}
\item Hom$(A)$ (\textit{homogeneity holds over} A) if  whenever $B \cup C \subseteq \A$ is disjoint from $A$ and $B \sim C$, there is an automorphism $\pi$ of $U$ such that $\pi B = C$ and $\pi$ pointwise fixes $A$.

\item A set $B$ of urelements \textit{duplicates} $A$ if $B$ and $A$ are disjoint and $A \sim B$.

\item Dup$(A)$  (\textit{duplication holds over} A) if  whenever a set $B$ of urelements is disjoint from $A$, there is a duplicate $C$ of $B$ that is disjoint from $A$.

\item Dup$^+(A)$  (\textit{strong duplication holds over} A) if whenever a set $B$ of urelements is disjoint from $A$, there is an infinite family $F$ of pairwise disjoint sets of urelements such that each $C \in F$ is disjoint from $A$ and has a subset that duplicates $B$.
\end{enumerate}
\end{definition}
\noindent Intuitively, Hom$(A)$ says that all equinumerous sets of urelements outside $A$ are indistinguishable, and we will often use Dup$(A)$ to get Hom$(B)$ for every $B \supseteq A$ by the next proposition. The notion of strong duplication Dup$^+(A)$ is isolated by Elliot Glazer, and it will be used to prove our key lemma.
\begin{prop}\label{prop:dup(A)->hom(A)}
Let $A$ and $B$ be sets of urelements.
\begin{enumerate}
\item Dup($A$) implies Hom(A).
\item If Dup($A$) and $A \subseteq B$, then Dup($B$). 
\item If Dup$^+(A)$ and $A \subseteq B$, then Dup$^+(B)$. 
\end{enumerate}
\end{prop}
\begin{proof}
(1) Let $B$ and $C$ be such that $B \sim C$ and $(B \cup C) \cap A = \emptyset$. By Dup($A$), there is some duplicate $D$ of $B \cup C$ that is disjoint from $A$. Let $\pi_1$ be an automorphism that only swaps $B$ and some subset of $D$, and $\pi_2$ only swap $\pi_1 B$ and $C$. Then $\pi_2 \circ \pi_1$ is as desired.

(2) If $C$ is disjoint from $B$, then by duplicating $(B - A) \cup C$ over $A$ we can obtain a duplicate of $C$ that is disjoint from $B$.

(3) If $C$ is disjoint from $B$, then there is a family $F$ of pairwise disjoint sets of urelements such that each $D \in F$ is disjoint from $A$ and has a subset duplicates $(B - A) \cup C$. Then $F^* = \{D - B \mid D \in F\}$ is the desired family of sets.  \end{proof}

\begin{lemma}\label{lemma:dup(A)basics}
Let $A$ be a set of urelements. Suppose that every set of urelements disjoint from $A$ is well-orderable. Then
\begin{enumerate}
\item Dup$(C)$ for some $C \supseteq A$. 
\item If Collection holds, then Dup$^+(C)$ for some $C \supseteq A$.
\end{enumerate}
\end{lemma}
\begin{proof}
(1) For any $B \supseteq A$, an initial ordinal $\kappa$ is the \textit{tail cardinal} of $B$ if (i) $\kappa$ is realized by some $B'$ that is disjoint from $B$ and (ii) every $C$ disjoint from $B$ injects into $\kappa$.

Case 1: Every $B \supseteq A$ has a tail cardinal. Let $C \supseteq A$ be such that $C$ has the least tail cardinal $\lambda$ and we show that Dup($C$). If $D$ is disjoint from $C$, then there is some $E$ disjoint from $C \cup D$ which realizes the tail cardinal of $C \cup D$. Then $D \preceq \lambda \preceq E$ so $E$ contains a duplicate of $D$. 

Case 2: Some $C \supseteq A$ does not have the tail cardinal. We show that Dup($C$). Suppose that $B$ is infinite and disjoint from $C$. Then there is some $D$ disjoint from $C$ that does not inject into $B$. Then $D - B$ must contain a duplicate of $B$ since $B$ and $D$ are well-orderable.

(2) We may assume that $\A$ is a proper class. Then by Collection, for every $C \subseteq \A$ there is a well-orderable infinite $D \subseteq \A$ that is disjoint from $C$. By (1), Dup$(C)$ for some $C \supseteq A$. Suppose $B$ is disjoint from $C$. There is some well-orderable infinite $D$ disjoint from $C \cup B$ which contains a subset that duplicates $B$. Since $D \sim \omega \times D$,  it follows that Dup$^+(C)$.\end{proof}
\noindent It is consistent with $\ZFU$ that duplication holds over no set of urelements (see the model in \cite[Theorem 10]{YaoAbstraction}) and consistent with $\ZFCU$ that strong duplication holds over no set of urelements (see the model in \cite[Theorem 2.17 (3)]{YaoAxiomandForcing}).

\begin{lemma}\label{lemma:Plenitude+->Dup}
Assume Plenitude$^+$. Then for every set of urelements $A$,  Dup$^+(A)$.
\end{lemma}
\begin{proof}
We show that Dup$^+(\emptyset)$. Let $A$ be a set of urelements and $\kappa =\max \{\omega, \aleph(P(A))\}$. By Plenitude$^+$, there is a bijection $f$ from $A \times \kappa$ to some set $E$ of urelements. For each $\alpha < \kappa$, let $E_\alpha = f[A\times\{\alpha\}]$. The set $\{\alpha < \kappa \mid E_\alpha \cap A \neq \emptyset \}$ must have size less than $\kappa$, since otherwise the map $\alpha \mapsto E_\alpha \cap A$ will induce an injection from $\kappa$ to $P(A)$. So $\{\alpha < \kappa \mid E_\alpha \cap A = \emptyset \}$ has size $\kappa$ and it follows that Dup$^+(\emptyset)$.
\end{proof}

\begin{definition}
Let $\varphi$ be a formula with subformulas $\varphi_1, ..., \varphi_n$. For sets $u$ and $v$, $v$ is a $\varphi$-cover of $u$ if and only if (i) $v$ is transitive; (ii) $u \subseteq v$; and (iii) for every subformula $\varphi_i (x_1, ..., x_{m+1})$ of $\varphi$,
$$\forall x_1, ..., x_m \in u [\exists y \varphi_i (x_1, ..., x_m, y)  \to \exists y \in v \varphi_i (x_1, ..., x_m, y)].$$
\end{definition}

\begin{prop}\label{prop:Collection->Phicover}
Assume Collection. Let $\varphi$ be a formula. Then every set has a $\varphi$-cover.
\end{prop}
\begin{proof}
Let $u$ be a set and $\varphi_1, ..., \varphi_n$ be the subformulas of $\varphi$. For each subformula  $\varphi_i (x_1, ..., x_{m+1})$ of $\varphi$,  define $$u_i = \{s \in u^m \mid \exists y \exists x_1, ..., x_m \in u (\varphi_i (x_1, ..., x_m, y) \land s = \langle x_1, ..., x_{m} \rangle)\}.$$ 
By Collection, for each $i$ there is a $v_i$ such that for every $s \in u_i$, there is a $y \in v_i$ such that  $\exists x_1, ..., x_m \in u [\varphi_i (x_1, ..., x_m, y) \land s = \langle x_1, ..., x_m \rangle]$. Set $v = TC(\bigcup_{i\leq n} v_i \cup x)$. $v$ is then a $\varphi$-cover of $u$. \end{proof}

\begin{lemma}\label{lemma:HP+Dup+Col->RP}
Suppose that cardinality is parametrically definable and Dup$^+(A)$ for some set $A$ of urelements. Then Collection implies RP.\footnote{Stronger versions of this lemma appear in \cite[Theorem 12]{HamkinsYaoRP2} and in an earlier draft of this paper. Both proofs contain errors. The current proof fixes these errors at the end by using strong duplication, and this modification is due to Elliot Glazer.}
\end{lemma}
\begin{proof}
By Lemma \ref{lemma:HP->AlmostUniversalSet} and \ref{lemma:dup(A)basics}, it follows that there is a set $A$ of urelements such that Dup$^+(A)$ and $V(A)$ contains every cardinal. Moreover, we can always find an arbitrarily tall $V_\gamma(A)$ that is closed under cardinal addition.

\begin{claim}
for every ordinal $\alpha$, there is an ordinal $\gamma \geq \alpha$ such that for every sets $x$ and $y$, $|x \cup y| \in V_\gamma (A)$ whenever $|x|$ and $|y|$ are in $V_\gamma (A)$.
\end{claim}
\begin{claimproof}
For every ordinal $\alpha$, $\bigcup_{\carda, \cardb \in V_\alpha(A)} X_{\carda, \cardb}$ is a set by Lemma \ref{lemma:SVC->UniSetandHP} and hence a subset of some $V_\beta (A)$. This allows us to define a sequence of ordinals $\<\gamma_n \mid n < \omega>$ by recursion such that $\gamma_0 = \alpha$ and $\gamma_{n+1}$ is the least $\beta$ such that $X_{\carda, \cardb} \subseteq V_{\beta}(A)$ for every $\carda$ and $\cardb$ in $V_{\gamma_n}(A)$. Now let $\gamma = \sup_{n <\omega} \gamma_n$. For every $|x|, |y| \in V_{\gamma_n}(A)$, $|x \cup y| \in X_{|x|, |y|} \subseteq V_{\gamma_{n+1}}(A)$.\end{claimproof}
\vspace{4pt}

 Now assume Collection and fix some formula $\varphi$. We show that RP holds for $\varphi$. Fix some set $x$. We may assume that $x$ is in some $V_{\lambda_0}(A)$, as we can enlarge $A$ if needed. We first construct an ordinal that is tall enough, which will be the height of the transitive set that reflects $\varphi$

$\carda$ is a \textit{disjoint cardinal} (for $A$) if there is some $B \subseteq \A$ such that $|B| = \carda$ and $B \cap A = \emptyset$ (we say that $B$ witnesses $\carda$).
\begin{claim}\label{claim:DisjointCardinal}
Let $\alpha$ be an ordinal. For every disjoint cardinal $\carda$ in $V_\alpha(A)$, there is some ordinal $\beta$ with the following property.
\begin{enumerate}
\item [] (*) There is a disjoint cardinal $\cardb \in V_\beta(A)$ such that for every $B$ that witnesses $\carda$, there is some $B^+ \supseteq B$ witnessing $\cardb$ such that $V_\beta (A \cup B^+)$ is a $\varphi\text{-cover of }  V_\alpha(A \cup B)$.
\end{enumerate}
\end{claim}
\begin{claimproof}
Let $B'$ witness $\carda$. By Proposition \ref{prop:Collection->Phicover}, $V_\alpha(A \cup B')$ has a $\varphi$-cover, which is in $V_\gamma (A \cup B^+)$ for some ordinal $\gamma$ and $B^+ \supseteq B'$ that is disjoint from $A$. Let $\cardb = |B^+|$, which is in $V_\beta (A)$ for some $\beta \geq \gamma$. So $V_\beta (A \cup B^+)$ a $\varphi$-cover of $V_\alpha(A \cup B')$.  Now if $B$ is any witness of $\carda$, by Hom(A) there is an automorphism $\pi$ such that $\pi B' = B$ and $\pi$ pointwise fixes $A$. It follows that $V_\beta (A \cup \pi B^+)$ is a $\varphi$-cover of $V_\alpha(A \cup B)$, $B \subseteq \pi B^+$ and $|\pi B^+| = \cardb$.
\end{claimproof}

Thus, for every ordinal $\alpha$ and disjoint cardinal $\carda \in V_\alpha (A)$, let $\beta_\carda$ be the least ordinal such that
\begin{enumerate}
\item $\beta_\carda$ has the property (*) as in  Claim \ref{claim:DisjointCardinal}; and

\item $|x \cup y| \in V_{\beta_\carda} (A)$ whenever $|x|$ and $|y|$ are in  $V_{\beta_\carda} (A)$.
\end{enumerate}
Then we define an increasing sequence of ordinals $\<\lambda_n \mid n < \omega>$ such that $\lambda_{n+1} = \sup \{\beta_\carda \mid \carda \text{ is a disjoint cardinal in } V_{\lambda_n}(A) \}$. Let $\lambda = \sup_{n<\omega} \lambda_n$. 

We next go sideways to construct a set of urelements that is big enough. By applying Collection to the set of disjoint cardinals in $V_\lambda(A)$, we get a set $C$ of urelements such that $C$ is disjoint from $A$ and $\carda \leq |C|$ for every disjoint cardinal $\carda \in V_\lambda(A)$. By Dup$^+(A)$, there is an infinite family $F$ of pairwise disjoint sets of urelements such that for every $B \in F$, $B$ is disjoint from $A$ and $B$ has a subset that duplicates $C$. Define
$$Z = \{ B \subseteq \bigcup F \mid \exists S \in [F]^{<\omega} (B \subseteq \bigcup S \land |B| \in V_\lambda(A))\}.$$
$Z$ is closed under finite union. Now let $t = \bigcup_{B \in Z} V_\lambda (A \cup B)$.  We show that $t$ reflects every subformula of $\varphi$ by induction on $\varphi_i$. 

The key step is to verify that $t$ passes the Tarski-Vaught test. So suppose that $x_1, ..., x_m \in t$ and $\exists y \varphi_i (x_1, ..., x_m, y)$. We need to show that $\exists y \in t$ with $\varphi_i (x_1, ..., x_m, y)$. $x_1, ... x_m$ are in some $V_{\lambda_n} (A \cup B)$ for some $\lambda_n$ and $B \in Z$, where $|B| \in V_{\lambda_n}(A)$ and $B \subseteq \bigcup S$ for some finite subset $S$ of $F$. So there is some disjoint cardinal $\cardb \in V_{\lambda_{n+ 1}} (A)$ witnessed by some $B^+ \supseteq B$ such that $V_{\lambda_{n+ 1}} (A \cup B^+)$ is a $\varphi$-cover of $V_{\lambda_n} (A \cup B)$, and hence there is some $y \in V_{\lambda_{n+ 1}} (A \cup B^+)$ with $\varphi_i (x_1, ..., x_m, y)$. It follows that there is some $B^* \in F$ such that $B^* \notin S$ and hence $B^* \cap (A \cup B) = \emptyset$. Since $B^*$ has a subset duplicating $C$, $\cardb \leq |B^*|$. By Hom($A \cup B$), there is an automorphism $\pi$ such that $\pi B^+  \subseteq B \cup B^*$ and $\pi$ pointwise fixes $A \cup B$. $\pi B^+ \subseteq \bigcup (S \cup \{B^*\})$ and $|\pi B^+| = \cardb$ so $\pi B^+ \in Z$. Thus, $\varphi_i(x_1, ..., x_m, \pi y)$ and $\pi y \in t$. This completes the proof.\end{proof}

\noindent Lemma \ref{lemma:Plenitude+->Dup} and  \ref{lemma:HP+Dup+Col->RP} then yield Theorem 1.
\begin{manualtheorem}{1}\label{thm:HP+Ple+Col->RP}
Assume that cardinality is parametrically definable and Plenitude$^+$ holds. Collection is equivalent to RP. \qed
\end{manualtheorem}

\subsection{Universal set, cardinal representative, and SVC}
We now consider some choice principles under which Plenitude$^+$ yields RP.
\begin{definition}
A set of urelements $A$ is \textit{universal} if every set $x$ injects into $V(A)$.
\end{definition}
\noindent The existence of a universal set is also a height-maximality principle: it says that the universe is built tall enough from some set of urelements. It is also a choice principle: AC implies that $\emptyset$ is universal, while it is consistent with $\ZFU$ that no universal set exists (see the model in \cite[Theorem 10]{YaoAbstraction}).

\begin{prop}\label{prop:EquivsofUni(A)}
Let $A$ be a set of urelements. The following are equivalent.
\begin{enumerate}
\item $A$ is universal.
\item Every set of urelements disjoint from $A$ injects into $V(A)$.
\item Every set of urelements injects into $V(A)$. 
\end{enumerate}
\end{prop}
\begin{proof}
(1) $\to$ (2) is clear.

(2) $\to$ (3). Assume (2) and let $B$ be a set of urelements. Since $B - A$ injects into $V(A)$, $B - A$ is equinumerous to $\{0\} \times x$ for some $x \in V(A)$. Therefore, $B$ injects into $A \cup (\{0\} \times x)$, which is a set in $V(A)$.

(3) $\to$ (1). Let $x$ be a set. By (2), there is an injection $g$ from $ker(x)$ into $V(A)$. Define a function $f$ on $TC(x)$ by recursion as follows.
\begin{equation*}
f(x) = 
\begin{cases*}
g(x)   & $x \in ker(x)$ \\
\{f(y) \mid y \in x\} & otherwise
\end{cases*}
\end{equation*}
A straightforward induction shows that $f$ is an injection from $TC(x)$ into $V(A)$.
\end{proof}

\begin{prop}\label{prop:UniSet->HP}
If there is a universal set of urelements, then cardinality is parametrically definable.
\end{prop}
\begin{proof}
Let $A$ be universal and we simply use Scott's trick within $V(A)$. That is, for every set $x$, let $$|x| = \{y \in V(A) : y \sim x \land \forall z \in V(A) (z \sim x \rightarrow \rho (y) \leq \rho(z))\},$$
where $\rho (z)$ is the least ordinal $\alpha$ such that $z \in V_\alpha (A)$.
$|x|$ is a set and it is routine to check that the map $x \mapsto |x|$ is a cardinality assignment function.
\end{proof}
Glazer \cite{GlazerMOFUnboundedCards} proved that the definability of cardinality does not imply the existence of a universal set.

\begin{corollary}\
Assume that every set of urelements is well-orderable. Then Collection is equivalent RP.
\end{corollary}
\begin{proof}
Under the assumption, every set of urelements injects into V, so $\emptyset$ is universal and hence cardinality is definable. If Collection holds, we will have Dup$^+(\emptyset)$ by Lemma \ref{lemma:dup(A)basics} (2) and so RP holds by Lemma \ref{lemma:HP+Dup+Col->RP}.
\end{proof}

\begin{lemma}\label{lemma:Uni(A)+Ple+->RP}
If there is a universal set of urelements, then Plenitude$^+$ implies RP.
\end{lemma}
\begin{proof}
Suppose that there is a universal set $A$ of urelements and Plenitude$^+$ holds. Then cardinality is parametrically definable and Dup$^+(B)$ holds for every set of urelements $B$ by Lemma \ref{lemma:Plenitude+->Dup}. So by Lemma \ref{lemma:HP+Dup+Col->RP} it remains to prove Collection.

Assume that $\forall x \in w \exists y \varphi(x, y, p)$ for some set $w$ and object $p$. Let $B = ker (w) \cup ker(p)$. For every $x \in w$, let $\alpha_x$ be the least ordinal such that $V_{\alpha_x} (A)$ contains some $z$ such that $z \sim ker(y)$ for some $y$ with $\varphi(x, y, p)$. Let $\gamma = \sup_{x \in w} \alpha_x$. There is then some $C' \subseteq \A$ that realizes $V_\gamma(A)$. By Dup($\emptyset$), we can duplicate $B \cup C'$ to get a $C$ such that $C \sim C'$ and $C \cap B = \emptyset$. By Proposition \ref{prop:CollectionEquivs} it remains to show that $\forall x \in w \exists y \in V(B \cup C) \varphi(x, y, p)$. Let $x \in w$. Then there is some $z \in V_\gamma(A)$ such that $z \sim ker(y)$ for some $y$ with $\varphi(x, y, p)$. It follows that $ker(y) - B$ injects into $C$. So by Hom($B$), there is an automorphism $\pi$ that moves $ker(y)$ into $B \cup C$ while pointwise fixing $B$. Thus, we have $\varphi(x, \pi y, p)$ and $ker(\pi y) \subseteq B \cup C$. This completes the proof.\end{proof}

A natural strengthening of the definability of cardinality is to ask the cardinality assignment function also to select a \textit{representative} from the class $\{y \mid y \sim x\}$ for each set $x$.
\begin{definition}
Cardinality is (parametrically) \textit{representable} if there is a (parametrically) definable cardinality assignment function $x \mapsto |x|$ such that $x \sim |x|$ for every set $x$.\end{definition}
\noindent AC implies that cardinality is representable. Pincus \cite{pincus1974cardinal} proves that ZF cannot prove cardinality is representable and that ZF + ``cardinality is representable'' does not prove AC.

\begin{manualtheorem}{2}\label{thm:HP+andPle+->RP}
Assume that cardinality is parametrically representable. Plenitude$^+$ implies RP.
\end{manualtheorem}
\begin{proof}
Suppose that cardinality is parametrically representable. By Lemma \ref{lemma:HP->AlmostUniversalSet}, there is a set $A$ of urelements with $|x| \in V(A)$ for every set $x$. But $x \sim |x|$ for every set $x$ so every set injects into $V(A)$, making $A$ a universal set. The theorem then follows from Lemma \ref{lemma:Uni(A)+Ple+->RP}.
\end{proof}

Finally, consider the principle of \textit{Small Violations of Choice} introduce by Blass \cite{blassInjec1979}, which says that AC is only violated in some local manner.
\begin{definition}
Let $S$ be a set. SVC(S) holds if for every set $x$ there is an ordinal $\alpha$ such that $x \preceq^* S \times \alpha$. SVC abbreviates $\exists S\ $SVC$(S)$.
\end{definition}
\noindent SVC has been widely studied in the literature, and it has numerous interesting consequences in ZF. For example, as shown by Blass \cite[Theorem 4.6]{blassInjec1979}, SVC is equivalent to the statement that AC holds in a forcing extension (and this equivalence holds in $\ZFU$).

\begin{lemma}\label{lemma:SVC->UniSetandHP}
SVC implies the existence of a universal set of urelements.
\end{lemma}
\begin{proof}
Suppose that SVC(S) holds for some set $S$. We show that $ker(S)$ is universal.  Let $x$ be any set. Then there is some ordinal $\alpha$ such that $S \times \alpha$ maps onto $x$. So $x$ maps injectively into $P(S \times \alpha)$ which is in $V(ker(S))$. \end{proof}
\noindent Our last positive result immediately follows from Lemma \ref{lemma:Uni(A)+Ple+->RP} and \ref{lemma:SVC->UniSetandHP}.
\begin{manualtheorem}{3}\label{thm:SVC+Ple+->RP}
Assume SVC. Plenitude$^+$ implies RP. \qed
\end{manualtheorem}
\noindent We conjecture that over $\ZFU$, SVC $+$ Collection implies RP. In \cite[Theorem 10]{YaoAbstraction}, the author proves that given any initial ordinal $\kappa$, it is consistent with $\ZFU$ + DC$_\kappa$-scheme + RP that cardinality is not parametrically definable. Thus, SVC is not provable from  $\ZFU$ + DC$_\kappa$-scheme + RP for any given initial ordinal $\kappa$.

\section{Independence Results}\label{section:IndependenceResults}

\subsection{Permutation and small-kernel models} 
The independence results will be proved by using inner models of permutation models. We begin by reviewing the basics of permutation models. As we will start with models of $\ZFCU$ where the urelements might be a proper class, some verifications are needed. However, the basic ideas remain the same, and we refer the reader to \cite{jech2008axiom} for a classic exposition.

\begin{definition}\label{permutationmodeldef}
Let $A$ be a set of urelements and $\G_A$ be a \textit{group} of permutations of $A$. For every $x$, define $sym(x) = \{\pi \in \G_A \mid \pi x = x\}$; if $x$ is a set, define $fix(x) = \{\pi \in \G_A \mid \pi y = y \text{ for all } y\in x \}$. A \textit{normal filter} $\F$ on $\G_A$ is a nonempty set of subgroups of $\G_A$ which contains $sym(a)$ for every urelement $a$ and is closed under supergroup, finite intersection, and conjugation (i.e., for all $\pi \in \G_A$ and $H \in \F$, $\pi H \pi^{-1} \in \F$). An object $x$ is \textit{symmetric} (with respect to $\F$) if $sym(x) \in \F$. The \textit{permutation model} $W$, determined by $\F$, is the class of all \textit{hereditarily symmetric} objects, i.e., $W = \{x \mid x \text{ is symmetric} \land x \subseteq W\}$.
\end{definition}
\noindent A standard argument shows that for every $\pi \in \G_A$ and $x \in W$, $\pi x \in W$. So every $\pi \in \G_A$ is an automorphism of $W$. Moreover, since the notion of symmetry only concerns the set of urelements $A$, every object $x$ whose kernel is disjoint from $A$ is also in $W$ since $fix(x) = \G_A$. In particular, if $B$ is a well-orderable set of urelements disjoint from $A$, it remains well-orderable in $W$.

We now review the fundamental theorem of permutation models and verify that it preserves the additional axioms discussed previously.
\begin{theorem}\label{thm:FundamentalThmofPermutationModel}
Let $A$, $\G_A$ and $\F$ be as in Definition \ref{permutationmodeldef} and $W$ be the resultant permutation model. Then
\begin{enumerate}
    \item $W \models \ZFU$;
    \item $W \models$ Collection if $U \models$ Collection.
    \item $W \models$ Plenitude if $U \models$ Plenitude.
    \item $W \models $ RP if $U \models$ RP + AC.
\end{enumerate}
\end{theorem}
\begin{proof}
(1) Since $W$ is transitive and contains all pure sets, Extensionality, Foundation, and Infinity all hold in $W$. Union holds in $W$ because for any set $x \in W$, $sym(x) \subseteq sym(\bigcup x)$. If $x, y \in W$, $sym(x) \cap sym(y) \subseteq sym(\{x, y\})$ and so $W$ satisfies Pairing. To show that $W$ satisfies Powerset, let $x \in W$ be a set. One observes that $P^W(x) = \{y \in W : y \subseteq x\}$ is symmetric since $sym(x) \subseteq sym(P^W(x))$.

\textit{Separation.} Let $x \in W$. We show that $v = \{y \in W : y \in x \land \varphi^W(y, u)\}$ is symmetric, where $u$ is a parameter in $W$. If $\pi \in sym(x) \cap sym(u)$ and $y \in v$, it follows that $\pi y \in v$ since $\varphi^W(\pi y,x,u)$. So $sym(x) \cap sym(v) \subseteq sym(v)$ and hence $v$ is symmetric.

\textit{Replacement.} Suppose that $W \models \forall x \in w \exists ! y \varphi(x, y, u)$, where $w, u \in W$. Let $v = \{y \in W : \exists x \in w \ \varphi^W(x, y, u)\}$ and we show that $sym(w)\cap sym(u) \subseteq sym(v)$. If $\pi \in sym(w)\cap sym(u)$ and $y \in v$, then $\varphi^W(x, y, u)$ for some $x \in w$ and so $\varphi^W(\pi x, \pi y, u)$ for some $\pi x \in w$; thus, $\pi y \in v$ and hence $sym(w)\cap sym(u) \subseteq sym(v)$.

(2) Suppose that Collection holds in $U$ and $W \models \forall x \in w \exists y \varphi(x, y, u)$ for some $w, u \in W$. So in $U$ there is a $v$ such that $\forall x \in w \exists y \in v (y \in W \land \varphi^W(x, y, u))$. Define
$$v' = \{\pi y \mid y \in W \cap v \land \pi \in \G_A \}.$$
$v' \in W$ because $sym(v') = \G$. $v \subseteq v'$ so $v'$ is a collection set in $W$.

(3) Let $\kappa$ be an infinite initial ordinal. In $U$, let $B$ realize max$\{\aleph(A), \kappa\}$. Since $A \cap B \prec \aleph (A)$, $B - A$ contains a subset $C \in W$ such that $(C \sim \kappa)^W$.

(4) Since $U \models$ Collection, $W \models $ Collection by (2). By Lemma \ref{lemma:HP+Dup+Col->RP} and Proposition \ref{prop:UniSet->HP}, it remains to show that in $W$, there is a universal set $C$ such that Dup$^+(C)$. Because every set of urelements disjoint from $A$ is well-orderable in $W$, every $C \supseteq A$ is universal in $W$; and for some $C \supseteq A$ we have $W \models$ Dup$^+(C)$ by Lemma \ref{lemma:dup(A)basics} (2). \end{proof}
\noindent The assumption that $U \models$ AC in (4) should not be necessary but this weaker result will suffice for our purpose. As an easy application, we use permutation models to separate Plenitude and Plenitude$^+$.

\begin{theorem}
There is a model of $\ZFU$ in which
\begin{enumerate}
\item Plenitude holds;

\item RP holds;

\item Dup($\emptyset$) fails.
\end{enumerate}
Therefore, $\ZFU$ + Plenitude + RP does not prove Plenitude$^+$.
\end{theorem}
\begin{proof}
Let $U$ be a model of $\ZFCU$ + Plenitude. In $U$, fix some infinite set $A$ of urelements and let $\G_A$ be the group of all permutations of $A$. Let $\F$ be the normal filter generated by the finite ideal of $A$, i.e.,
$$\F = \{ H \subseteq \G_A: H \text{ is a subgroup of }\G_A \text{ and } fix(E) \subseteq H \text{ for some finite } E \subseteq A \}.$$
In the resultant permutation model $W$, both Plenitude and RP hold. Dup($\emptyset$) fails in $W$: a standard argument shows that $A$ is not well-orderable as $fix(A) \notin \F$ but every set of urelements disjoint from $A$ is. So it follows from Lemma \ref{lemma:Plenitude+->Dup} that Plenitude$^+$ fails in $W$. \end{proof}

However, Theorem \ref{thm:FundamentalThmofPermutationModel} tells us that to get a failure of Collection we need to create a somehow deeper failure of AC. The idea here is to create a proper class of non-well-orderable sets of urelements using $A$.
\begin{definition}
A \textit{class} $\I$ of sets of urelements is an $\A$-\textit{ideal} if it is a class ideal on $\A$ (which does not contain $\A$ if it is a set) that includes all singletons of urelements. Given an $\A$-ideal $\I$, the \textit{small kernel model} $U^\I$ is the class $\{x \in U : ker(x) \in \I\}$.
\end{definition}

\begin{theorem}
Assume $\ZFU$. For every $\A$-\textit{ideal} $\I$, $U^\I \models \ZFU$.
\end{theorem}
\begin{proof}
See \cite[Theorem 2.16]{YaoAxiomandForcing}, where the assumption that $U \models$ AC is clearly not needed for the preservation of other axioms.
\end{proof}
\noindent Small-kernel models have appeared in many places including \cite{levy1969definability}, \cite{blassInjec1979}, \cite{FinitenessSVC} and \cite{YaoAbstraction}. Unlike permutation models, small-kernel models preserve neither Collection nor RP (\cite[Theorem 2.17]{YaoAxiomandForcing}). Moreover, it can be used to destroy SVC (\cite[Theorem 3.1]{blassInjec1979}) and the definability of cardinality (\cite{levy1969definability} and \cite{YaoAbstraction}).

\subsection{A general construction}
We generalize L\'evy's model in \cite[Section 4]{levy1969definability}, where he proves that in $\ZFU$ cardinality can be parametrically definable but not definable. Let $U$ be a model of $\ZFCU$ + Plenitude. Fix a set of urelements $A$ of size $\kappa$ for some infinite initial ordinal $\kappa$, and enumerate $A$ with $\kappa \times \Q$. $A$ is thus the matrix $\bigcup_{\alpha <\kappa} A_\alpha$ and for each row $A_\alpha$, there is a relation $<_\alpha$ such that $\<A_\alpha, <_\alpha>$ is isomorphic to the rationals $\<\Q, <_\Q >$. We write $A_\alpha = \{a^\alpha_i : i \in \Q\}$ for each $\alpha <\kappa$. 
A permutation $\pi$ of $A$ is said to be a \textit{uniform automorphism} of $A$ if there is an automorphism $\rho$ of $\<\Q, <_\Q>$ such that $\pi(a^\alpha_i) = a^\alpha_{\rho(i)}$ for every $\alpha < \kappa$ and $i \in \Q$, i.e., $\pi$ follows the same automorphism of $\Q$ in each row. Let $\G_A$ be the group of all uniform automorphisms of $A$. Define $\F$ to be the normal filter generated by the finite ideal of $A$, i.e., 
$$\F = \{ H \subseteq \G_A: H \text{ is a subgroup of }\G_A \text{ and } fix(E) \subseteq H \text{ for some finite } E \subseteq A \}.$$
Call the resultant permutation model $W$. Let $\I$ be the following $\A$-ideal.
$$\I = \{B \subseteq \A \mid \{\alpha < \kappa \mid B \cap A_\alpha \neq \emptyset \} \text{ has size} < \kappa \}.$$
The corresponding small-kernel model $W^\I$, which only contains objects whose kernel intersects $A$ in a small block of rows, is then a model of $\ZFU$. We now prove several general facts about $W^\I$.

$W^\I \models $ Plenitude since in $U$ every initial ordinal $\kappa$ is realized by some set of urelements disjoint from $A$. For any urelement $a^\alpha_i$, a function $f$ is said to \textit{vertically fix} $a^\alpha_i$ if $f (a^\alpha_i) = a^\beta_i$ for some $\beta < \kappa$.
\begin{lemma}\label{Anareduplicates}
\
\begin{enumerate}
\item For every $\alpha, \beta < \kappa$, $W^\I \models A_\alpha \sim A_\beta$.
\item $W^\I \models $ Dup$(\emptyset)$; and $W^\I \models $ Dup$^+(\emptyset)$ if $\kappa > \omega$.
\item For every $\alpha < \kappa$, $W^\I \models A_\alpha$ is universal.
\end{enumerate}
\end{lemma}
\begin{proof}
(1) For every $\alpha, \beta < \kappa$, in $U$ let $f_{\alpha, \beta}$ be the bijection between $A_\alpha$ and $A_\beta$ such that $f$ vertically fixes every urelement in $A_\alpha$.  Consider any $\pi \in \G_\A$ and let $\rho$ be the automorphism of $\Q$ followed by $\pi$. For every $\<a^\alpha_i, a^\beta_i> \in f$, $\<\pi a^\alpha_i, \pi a^\beta_i > = \<a^\alpha_{\rho i}, a^\beta_{\rho i}> \in f$. So $sym(f_{\alpha, \beta}) = \G_A$ and hence $f_{\alpha, \beta} \in W^\I$.

(2) For any initial ordinal $\lambda < \kappa$, any two $\lambda$-blocks of rows $\bigcup_{\eta<\lambda} A_{\alpha_\eta}$ and $\bigcup_{\eta<\lambda} A'_{\beta_\eta}$ are duplicates in $W^\I$ since the bijection $F = \bigcup_{n < \omega} f_{\alpha_\eta, \beta_\eta}$ is in $W^\I$. This shows that $W^\I \models$ Dup$(\emptyset)$.

Suppose $\omega < \kappa$. For any $\bigcup_{\eta<\lambda} A_{\alpha_\eta}$, where $\lambda < \kappa$, there is some $\bigcup_{\eta<\theta} A'_{\alpha_\eta}$ disjoint from  $\bigcup_{\eta<\lambda} A_{\alpha_\eta}$, where $\theta = \max \{\omega, \lambda\}$. So in $M$ we can partition $\bigcup_{\eta<\theta} A'_{\alpha_\eta}$ into $\omega$-many $\lambda$-block of rows, and each of them is a duplicate of $\bigcup_{\eta<\lambda} A_{\alpha_\eta}$ in $M$. Thus, $W^\I \models$ Dup$^+(\emptyset)$.

(3) By Proposition \ref{prop:EquivsofUni(A)}, it suffices to show that in $W^\I$ every set of urelements injects in to $V(A_\alpha)$. Let $B$ be a set of urelements in $W^\I$. $B - A$ is equinumerous to some ordinal $\beta$. $B \cap A$ is contained in $\bigcup_{\eta <\lambda} A_{\beta_\eta}$ for some $\lambda < \kappa$. As the sequence of bijections $\<f_{\beta_\eta, \alpha} \mid \eta < \lambda>$ is in $W^\I$, there is a bijection in $W^\I$ between $\bigcup_{\eta <\lambda} A_{\beta_\eta}$ and $\bigcup_{\eta <\lambda} (\{\eta\} \times A_\alpha)$, and the latter contains no pure set. Thus, $B$ injects into $\beta \cup \bigcup_{\eta <\lambda} (\{\eta\} \times A_\alpha)$, which is in $V(A_\alpha)$. \end{proof}

Although small-kernel models do not necessarily satisfy SVC, we will show that SVC holds in $W^\I$. The point is that Blass' argument (\cite[Theorem 4.6]{blassInjec1979}) still works within $V(A_\alpha)^{W^\I}$ but $A_\alpha$ is universal in $W^\I$.
\begin{lemma}\label{lemma:WmodelsSVC}
$W^\I \models $ SVC.
\end{lemma}
\begin{proof}
Fix any $A_\alpha$. we show that there is an $S \in V(A_\alpha)^W$ such that $$W^\I \models \forall x \in V(A_\alpha) \exists \beta (x \preceq^* S \times \beta).$$
It follows that SVC$(S)$ holds in $W^\I$ since $A_\alpha$ is universal in $W^\I$ by Lemma \ref{Anareduplicates}.

We simply follow Blass' argument. In $U$, define
$$\F' = \{H \in \F \mid \exists x \in V(A_\alpha)^W sym(x) = H \}.$$
By AC in $U$, for each $H \in \F'$ we can choose a pair $x_H = \<\beta_H, y_H>$ such that $y_H \in V(A_\alpha)^W$, $sym(y_H) = H$ and $\beta_H \neq \beta_{H'}$ whenever $H \neq H'$. As a result, for every $H \in \F'$, $x_H \in V(A_\alpha)^W$ and $\pi(x_H) \neq x_{H'}$ for any $\pi \in \G_A$ and any different two $H$ and $H'$ in $\F'$. Let $S = \{\pi (x_H) \mid \pi \in \G_A \text{ and } H \in \F'\}$. $S \in W$ and moreover, since $\pi A_\alpha = A_\alpha$ for each $\pi \in \G_A$, we have $S \in V(A_\alpha)^W$. Therefore, $S \in W^\I$. 

Consider any $x$ in $W^\I$ such that $x \in V(A_\alpha)$. In $U$, enumerate $x$ as $\{y_\eta \mid \eta < \beta\}$ for some $\beta$ and define
$$f = \{\<\pi (x_H), \eta, \pi (y_\eta)> \mid \pi \in \G_A, H \in \F', \eta < \beta \text{ and } sym(y_\eta) = H\}.$$
\noindent $f \subseteq W$ and $sym(f) = \G_A$, so $f \in W$. Furthermore, $f \in V(A_\alpha)$ as $x \in V(A_\alpha)$, so $f \in W^\I$. $dom(f) \subseteq S \times \beta$ so it remains to show that $f$ is a surjection onto $x$. 

To show $f$ is a function, suppose that $(\pi (x_H), \eta) = (\pi' (x_{H'}), \eta')$. Then $x_H = \pi^{-1} \circ \pi' (x_{H'})$ so by the construction of $x_H$, we have $H = H'$. So $\pi^{-1} \circ \pi' \in sym(x_H) = sym(y_\eta)$. Thus, $\pi (y_\eta) = \pi'(y_\eta) = \pi'(y_{\eta'})$. To show that $f$ is onto $x$, let $y \in x$, which is $y_\eta$ for some $\eta < \beta$. Then $sym(y) \in \F'$ and so $x_{sym(y)} \in S$. It follows that $\<x_{sym(y)}, \eta, y_\eta> \in f$ and hence $f(x_{sym(y)}, \eta) = y$.\end{proof}

\begin{lemma}\label{alignment}
For every $\alpha < \kappa$, every injection $f$ in $W$ such that $dom(f) \subseteq A_\alpha$ vertically fixes co-finitely many urelements in its domain.
\end{lemma}
\begin{proof}
Of course, we may assume that $f$ has an infinite domain. Let $E$ be a finite subset of $A$ such that $fix(E) \subseteq sym(f)$ and $I$ be the rational indexes appeared in $E$, i.e., $I = \{q \in \Q \mid a^\beta_q \in E \text{ for some } \beta < \kappa\}$. We show that $f$ vertically fixes every urelement in $dom(f) - \{a^\alpha_p \mid p \in I\}$, which is co-finite. Suppose \textit{for reductio} that $f(a^\alpha_i) = a^\beta_j$ and $i \neq j$ for some $a^\alpha_i \in dom(f) - \{a^\alpha_p \mid p \in I\}$. Let $\pi \in \G_A$ be an automorphism induced by some automorphism of $\Q$ that moves $i$ but pointwise fixes $I \cup \{j\}$. It follows that $\pi \in fix(E \cup \{a^\beta_j\})$ and so $f(a^\alpha_i) = f (a^\alpha_{\pi(i)})$, which contradicts the assumption that $f$ is injective.\end{proof}

\begin{lemma}\label{Anhasnoduplicates}
For every $\alpha < \kappa$, if $B$ and $B'$ are subsets of $A_\alpha$ in $W$ such that $B  \triangle B'$ is infinite, then $B$ and $B'$ are not equinumerous. Thus, no infinite subset of $A_\alpha$ is well-orderable in $W$.
\end{lemma}
\begin{proof}
Suppose, without loss of generality, that $B - B'$ is infinite. Then no injection in $W$ can move the urelements in $B - B'$ into $B'$. For, such an injection would change the rational indexes for infinitely many urelements in $A_\alpha$, which is impossible by Lemma \ref{alignment}.\end{proof}

\subsection{Independence results}
Now we prove Theorem \ref{thm:Independence1} using the construction just introduced.
\begin{theorem}\label{thm:Ple+Dup+SVC/->Collection}
There is a model of $\ZFU$ in which 
\begin{enumerate}
    \item Plenitude holds;
    \item Dup$(\emptyset)$ holds;
    \item SVC holds; and
    \item Collection fails.
\end{enumerate}
\end{theorem}
\begin{proof}
Let $U$ be a model of $\ZFCU$ + Plenitude. Fix a set of urelements $A$ of size $\omega$ and enumerate it with $\omega \times \Q$, i.e., $A = \bigcup_{n<\omega}A_n$ and each row $A_n$ is identified with $\Q$. Let $W$ be the permutation model generated by (i) the group of all uniform automorphisms of $A$ and (ii) the normal filter induced by the finite ideal of $A$. Furthermore, define
$$\I = \{B \subseteq \A \mid \{n < \omega \mid B \cap A_n \neq \emptyset \} \text{ is finite}\}.$$
The lemmas in the previous subsection all apply to the resultant small-kernel model $W^\I$. So it remains to show that Collection fails in $W^\I$. 

A set of urelements $B$ is said to have a \textit{nice n-partition} if there are $B_1, ..., B_n$ such that 
\begin{itemize}
    \item [] (i) $B = B_1 \cup ... \cup B_n$;
    \item [] (ii) for each two $k, l \leq n$, $B_k, B_l$ are non-well-orderable duplicates.
\end{itemize}
In $W^\I$, we have 
\begin{align}\label{1}
     \forall m < \omega \exists B \subseteq \A (B \text{ has a nice } m \text{-partition}).
\end{align}
This is because each pair of $A_n$ are duplicates by Lemma \ref{Anareduplicates} so the union of any $m$-block of $A_n$ will do.

\begin{lemma}\label{partition}
In $W^\I$, there is no set of urelements $B$ such that for every $n < \omega$, $B$ has a subset with a nice $n$-partition.
\end{lemma}
\begin{proof}
Suppose \textit{for reductio} that such a $B$ exists in $W^\I$. We may assume that $B \cap A \subseteq A_1, ..., A_n$ for some finite number $n$. Let $C \subseteq B$ be a set with a nice $n+1$-partition such that $C = C_1 \cup ... \cup C_{n+1}$, where for each two $k, l \leq n$, $C_k$ and $C_l$ are non-well-orderable duplicates. And for each $1 \leq k < n+1$, let $f_k \in W^\I$ be a bijection from $C_k$ to $C_{k+1}$. Define an $(n+1)$-sequence of pairs $\<m_{1}, D_1, >, ..., \<m_{n+1}, D_{n+1}>$ by recursion as follows.
\begin{itemize}
   \item [] $m_1$ is the least number $\leq n$ such that $C_1 \cap A_{m_1}$ is infinite; and
    \item [] $D_1 = C_1 \cap  A_{m_1}$.
\end{itemize}   
Then for each $k < n$,
 \begin{itemize}
    \item [] $m_{k+1}$ is the least natural number $\leq n$ such that $f_k[D_k] \cap A_{m_{k+1}}$ is infinite and $D_{k+1} =  f_k [\{ a \in D_k  \mid f_{k} \text{ vertically fixes } a\}] \cap A_{m_{k+1}}$.
\end{itemize}
To check this is well-defined, we show that for each $k \leq n+1$, $m_k$ exists, $D_k \subseteq A_{m_k} \cap C_k$ and $D_k$ is infinite by induction. Suppose that $D_k \subseteq A_{m_k} \cap C_k$ is infinite. Then by Lemma \ref{Anhasnoduplicates}, $f_k[D_k]$ is a non-well-orderable subset of $C_{k+1}$ so it must have an infinite intersection with $A_l$ for some $l \leq n$, which means $m_{k+1}$ exists. By Lemma \ref{alignment}, $f_k$ must vertically fix infinitely many $a \in D_k$ such that $f_k (a) \in A_{m_{k+1}}$. Thus, $D_{k+1}$ is infinite. 

To finish the proof, we show that for every $j < k \leq n+1$, $m_k \neq m_j$, which is a contradiction because for each $k \leq n+1$, $m_k \leq n$. Suppose \textit{for reductio} that $m_{i+1} = m_j$, where $j < i+1 \leq n$. Then there is a sequence of urelements $\<a_1, ..., a_{i+1}>$ such that for each $1 < l \leq i+1$,
\begin{itemize}
\item [] (i) $a_l \in D_l$;

\item [] (ii) $a_l = f_{l-1}(a_{l-1})$; and

\item [] (iii) $f_{l-1}$ vertically fixes $a_{l-1}$.

\end{itemize}
$a_l \in A_{m_l}$ for each $l$ so $a_{i+1}, a_j \in A_{m_j}$ by assumption; but they also have the same rational index and hence $a_{i+1} = a_j$. However, this means that the intersection of $C_j$ and $C_{i+1}$ is nonempty, which contradicts the assumption we started with. This proves the lemma. \end{proof}
Now suppose that
\begin{align}\label{2}
  W^\I \models \exists v \forall m < \omega \exists B \in v (B \subseteq \A \land B \text{\ has a nice } m\text{-partition}).
\end{align}
Then for every $m < \omega$, $ker(v)$ has a subset with a nice $m$-partition, contradicting Lemma \ref{partition}. Therefore, (\ref{1}) and $\neg$(\ref{2}) yield a failure of Collection in $W^\I$.\end{proof}
We note that Theorem \ref{thm:Independence2} can be proved using permutation models alone. But here we prove a stronger result by constructing a counter model in which Dup$^+(\emptyset)$ holds, and this model also shows that SVC does not imply the (parameter-free) definability of cardinality.
\begin{theorem}
There is a model of $\ZFU$ in which
\begin{enumerate}
\item Plenitude holds;

\item Dup$^+(\emptyset)$ holds;

\item SVC holds;

\item RP holds;

\item cardinality is not definable; and

\item Plenitude$^+$ fails.
\end{enumerate}

\end{theorem}
\begin{proof}
Let $U$ be a model of $\ZFCU$ + Plenitude. In $U$, fix a set of urelements $A$ of size $\omega_1$ and enumerate it with $\omega_1 \times \Q$, i.e., $A = \bigcup_{\alpha <\omega_1} A_\alpha$ and each row $A_\alpha$ is identified with $\Q$ such that $A_\alpha = \{a^\alpha_i : i \in \Q\}$. As before, let $W$ be the permutation model generated by (i) the group $\G_A$ of all uniform automorphisms of $A$ and (ii) the normal filter induced by the finite ideal of $A$. Furthermore, define
$$\I = \{B \subseteq \A \mid \{\alpha < \omega_1 \mid B \cap A_\alpha \neq \emptyset \} \text{ is countable}\}.$$
The lemmas in the previous subsection all apply to the resultant small-kernel model $W^\I$. So it remains to prove (4), (5) and (6).\\

\noindent (4) Since $A_\alpha$ is universal in $W^\I$ by Lemma \ref{Anareduplicates}, cardinality is parametrcially definable in $W^\I$. So by Lemma \ref{lemma:HP+Dup+Col->RP} it suffices to show that $W^\I \models$ Collection. Let $w \in W^\I$ and suppose that for every $\forall x \in w \exists y \varphi(x, y)$. By Proposition \ref{prop:CollectionEquivs}, it is enough to find some $B\subseteq \A$ in $W^\I$ such that $W^\I \models \forall x \in w \exists y \in V(B) \varphi(x, y)$. Let $B_1 = \bigcup_{n<\omega} A_{\alpha_n}$ be a countable block of rows that includes $ker(w)$ and fix some countable block $B_2 =\bigcup_{n<\omega} A_{\beta_n}$ that is disjoint from $B_1$, and let $B = B_1 \cup B_2$. Let $x \in w$ and $y' \in W^\I$ be such that $W^\I \models \varphi(x, y')$. $ker(y') - B_1$ is contained in some countable block. But any two countable blocks of rows are equinumerous in $W^\I$ since the union of the individual canonical bijections is symmetric (as in Lemma \ref{Anareduplicates}). Thus, $ker(y') - B_1$ injects into $B_2$. Since $W^\I \models$ Dup$(\emptyset)$, $W^\I \models$ Hom$(B_1)$ by Lemma \ref{lemma:dup(A)basics}. Thus, in $W^\I$ there is an automorphism that pointwise fixes $B_1$ (and hence $w$) and moves $y'$ into $V(B)$. Therefore, $W^\I \models \exists y \in V(B) \varphi(x, y)$.\\

\noindent (5) Suppose \textit{for reductio} that cardinality is definable in $W^\I$. Then by Lemma \ref{lemma:HP->AlmostUniversalSet} every cardinal in $W^\I$ is a pure set. Let $0 < i< j \in \Q$ and $\rho$ be an automorphism of $\Q$ such that $\rho (i) = j$ and $\rho (0) = 0$. Fix some $A_\alpha$ and let $B_{0, i} = \{a^\alpha_k \mid 0<k<i\}$ and $B_{0, j} = \{a^\alpha_k \mid 0<k<j\}$, which are in $W^\I$. Consider $\pi \in \G_A$ induced by $\rho$. It follows that $\pi B_{0, i} = B_{0, j}$. As $\pi$ is an automorphism of $W^\I$, we have 
$$ W^\I \models |B_{0, i}| = \pi | B_{0, i}| = |\pi  B_{0, i}|  = |B_{0, j}|$$
Thus, $W^\I \models B_{0, i} \sim B_{0, j}$, but this contradicts Lemma \ref{Anhasnoduplicates}.\\

\noindent (6) Fix an $\alpha < \omega_1$, we show that $A_\alpha \times \omega_1$ is not realized in $W^\I$.
\begin{claim}
In $W^\I$, if $Z \subseteq A_\alpha \times \omega_1$ is such that $ (A_\alpha \times \omega_1) - Z$ is well-orderable, then $\omega_1$ injects into $Z$.
\end{claim}
\begin{claimproof}
Note that $\bigcap_{\beta < \omega_1} \{ a \in A_\alpha \mid \<a, \beta> \in Z\}$ is nonempty. Otherwise, $a \mapsto \<a, \beta_a>$, where $\beta_a < \omega_1$ is the least ordinal such that $\<a, \beta_a> \notin Z$, is an injection from $A_\alpha$ into the complement of $Z$; but $A_\alpha$ is non-well-orderable in $W^\I$, which makes the complement of $Z$ non-well-orderable---contradiction. So fix some $a \in \bigcap_{\beta < \omega_1} \{ a \in A_\alpha \mid \<a, \beta> \in Z\}$. Then $\beta \mapsto \<a, \beta>$ is an injection from $\omega_1$ into $Z$.\end{claimproof}
\vspace{4pt}

\noindent Suppose \textit{for reductio} that $B$ is a set of urelements in $W^\I$ and $f \in W^\I$ is a bijection between $B$ and $A_\alpha \times \omega_1$. Since $B - A$ is well-orderable and $B$ is non-well-orderable, $B \cap A$ is non-well-orderable and contained in some countable block $\bigcup_{n<\omega} A_{\alpha_n}$. Then $f[B \cap A]$ is a subset of $A_\alpha \times \omega_1$ whose complement is well-orderable. So $\omega_1$ injects into $B\cap A$ and hence into $\bigcup_{n<\omega} A_{\alpha_n}$, which is a contradiction.\end{proof}

\section*{Funding}
The author was supported by NSFC No. 12401001 and the Fundamental Research Funds for the Central Universities, Peking University.
\printbibliography

@book{jech2008axiom,
  title={The axiom of choice},
  author={Jech, Thomas J},
  year={2008},
  publisher={Courier Corporation}
}

@incollection{levy1969definability,
  title={The definability of cardinal numbers},
  author={L{\'e}vy, Azriel},
  booktitle={Foundations of Mathematics},
  pages={15--38},
  year={1969},
  publisher={Springer}
}

@book{yao2023set,
  title={Set theory with urelements},
  author={Yao, Bokai},
  year={2023},
  eprint={2303.14274},
    archivePrefix={arXiv},
  publisher={University of Notre Dame}
}

@article{gauntt1967undefinability,
  title={Undefinability of cardinality},
  author={Gauntt, Robert J},
  journal={Lectures notes prepared in connection with the Summer Institute on Axiomatic Set Theory held at University of California, Los Angeles, IV-M},
  year={1967}
}

@article{pincus1974cardinal,
  title={Cardinal representatives},
  author={Pincus, David},
  journal={Israel Journal of Mathematics},
  volume={18},
  pages={321--344},
  year={1974},
  publisher={Springer}
}

@article{blassInjec1979,
 ISSN = {00029947},
 URL = {http://www.jstor.org/stable/1998165},
 author = {Andreas Blass},
 journal = {Transactions of the American Mathematical Society},
 pages = {31--59},
 publisher = {American Mathematical Society},
 title = {Injectivity, Projectivity, and the Axiom of Choice},
 urldate = {2025-08-18},
 volume = {255},
 year = {1979}
}

@article{YaoAxiomandForcing,
	author = {Bokai Yao},
	doi = {10.1017/jsl.2024.58},
	journal = {Journal of Symbolic Logic},
	title = {Axiomatization and Forcing in Set Theory with Urelements},
	year = {forthcoming}
}

@article{YaoAbstraction,
	author = {Bokai Yao},
        doi={10.1017/S1755020325100804},
	journal = {Review of Symbolic Logic},
	title = {Abstraction Principles and the Size of Reality},
	year = {forthcoming}
}

@article{FinitenessSVC,
	author = {Horst Herrlich and Paul Howard and Eleftherios Tachtsis},
	doi = {10.1215/00294527-3490101},
	journal = {Notre Dame Journal of Formal Logic},
	number = {3},
	pages = {375--388},
	title = {Finiteness Classes and Small Violations of Choice},
	volume = {57},
	year = {2016}
}

@article{HamkinsYaoRP2,
	author = {Joel David Hamkins and Bokai Yao},
	doi = {10.1017/jsl.2022.87},
	journal = {Journal of Symbolic Logic},
	number = {3},
	pages = {1007--1043},
	title = {Reflection in Second-Order Set Theory with Abundant Urelements Bi-Interprets a Supercompact Cardinal},
	volume = {89},
	year = {2024}
}

@MISC{GlazerYaoRPinZFU,
	author = {Elliot Glazer and Bokai Yao},
	title = {Reflection Principles in ZFU},
	HOWPUBLISHED = {Unpublished manuscript}
}

@MISC {GlazerMOFUnboundedCards,
    TITLE = {Is every set of cardinals bounded?},
    AUTHOR = {Elliot Glazer (https://mathoverflow.net/users/109573/elliot-glazer)},
    HOWPUBLISHED = {MathOverflow},
    NOTE = {URL:https://mathoverflow.net/q/501823 (version: 2025-10-20)},
    EPRINT = {https://mathoverflow.net/q/501823},
    URL = {https://mathoverflow.net/q/501823}
}

\end{document}